\documentclass[12pt]{article}

\usepackage[T1]{fontenc}
\usepackage{amssymb}
\usepackage{amsfonts}
\usepackage{amsmath}
\usepackage{mathrsfs}
\usepackage{graphicx}
\usepackage{amsbsy}
\usepackage{theorem}
\usepackage{color}
\usepackage{hyperref}
\usepackage[normalem]{ulem}
\newtheorem{theorem}{Theorem}[section]

\newtheorem{definition}[theorem]{Definition}

\newtheorem{remark}[theorem]{Remark}

\def\thetheorem{\thesection.\arabic{theorem}}
\def\thesection{\arabic{section}}

\def\theequation {\thesection.\arabic{equation}}
\def\beq{\begin{equation}\displaystyle}
\def\eeq{\end{equation}}
\def\bel{\begin{equation} \displaystyle \begin{array}{l} }
\def\eel{\end{array} \end{equation} }
\def\bell{\begin{equation} \displaystyle \begin{array}{ll}  }
\def\eell{\end{array} \end{equation} }

\def\bea{\begin{eqnarray}}
\def\eea{\end{eqnarray} }
\def\bean{\begin{eqnarray*}}
\def\eean{\end{eqnarray*} }

\catcode`@=11
\renewcommand\appendix{\bigskip {\noindent \Large \bf Appendix}
  \setcounter{section}{0}%
  \setcounter{subsection}{0}%
\setcounter{equation}{0}%
\setcounter{theorem}{0}%
\def\thetheorem{A.\arabic{theorem}}
\def\theequation {A.\arabic{equation}}}
\catcode`@=12

\def\RR{\mathbb{R}}

\def\bs{\bigskip}

\def\eps{\varepsilon}

\def\p{\partial}

\def\bchapo{\widehat{b}}

\begin{document}

\bs
\begin{center}
\Large \bf
A remark on duality solutions for some weakly nonlinear scalar conservation
laws
\end{center}

\bs

\begin{center}
F. James$^a$ and N. Vauchelet$^{b}$

\bs

 \medskip

 {\footnotesize $^a$ 
Universit\'e d'Orl\'eans, Math\'ematiques, Applications et Physique Math\'ematique d'Orl\'eans, \\
CNRS UMR 6628, MAPMO \\
F\'ed\'eration Denis Poisson, CNRS FR 2964, \\
45067 Orl\'eans Cedex 2, France}
 \\

 {\footnotesize $^b$ 
UPMC Univ Paris 06, UMR 7598, Laboratoire Jacques-Louis Lions, \\
CNRS, UMR 7598, Laboratoire Jacques-Louis Lions and \\
INRIA Paris-Rocquencourt, Equipe BANG \\
F-75005, Paris, France}
 \\

 \medskip
 {\footnotesize {\em E-mail addresses:} 
 {\tt Francois.James@univ-orleans.fr, vauchelet@ann.jussieu.fr}
 }

\end{center}

\bs

\begin{abstract}
We investigate existence and uniqueness of duality solutions
for a scalar conservation law with a nonlocal interaction kernel.
Following \cite{BJ99}, a notion of duality solution for such a 
nonlinear system is proposed, for which we do not have uniqueness.
Then we prove that a natural definition of the flux allows to
select a solution for which uniqueness holds.
\end{abstract}

\bs

\section{Introduction}
\label{intro}

At a continuous level, many physical or biological systems are modelled 
thanks to scalar conservation laws.
In this note we will focus on a weakly nonlinear system of the kind~:
\beq\label{eq:nl}
\p_t \rho + \p_x (a(u)\rho) = 0, \qquad \quad \p_x u = \rho, 
\eeq
where $a$ is a given smooth function, $a\in C^0(\RR)$.
This system is complemented with the initial data $\rho(t=0)=\rho^0$.
We notice that we can rewrite (\ref{eq:nl}) as a single equation
since we have $u=H*\rho$ where $H$ is the Heaviside function
and we recover the so-called non local aggregation equation.
This model arises in several applications in physics and biology
where a self-consistant interaction field $u$ governs the evolution of
a density of population $\rho$. Then $u$ is defined as
$u=-\p_x \phi$ where $\phi$ is the interaction potential.
For instance, in the modelling of cell movement by chemotaxis, $\phi$ is the 
concentration of some chemical called chemo-attractant (when $a$ is non-increasing)
or chemo-repellent (when $a$ is non-decreasing) which drives 
the dynamics of individuals (bacteria).
In gas dynamics, this model can be derived thanks to a high-field limit 
from the Vlasov--Poisson--Fokker--Planck system \cite{NPS}, a nonincreasing (resp. nondecreasing) $a$
corresponds the to the repulsive (resp. attractive) case.

From a mathematical viewpoint, it is well-known that in the attractive 
case, i.e. when $a$ is non-increasing, finite time blow-up of 
regular solutions for such system occurs
(see e.g. \cite{bertozzi} and references therein).
Therefore one has to look for solutions $\rho$ which are measure-valued in space, which generates
several difficulties, because the velocity $a(u)$ turns out to be discontinuous, so that
the product in the divergence term is not well defined, and the corresponding flow has to be defined cautiously.
A recent approach consists in using techniques from optimal transport, see \cite{Carrillo}. Another possibility
is to define {\it a priori} the product. For the Vlasov--Poisson system, this has been done in \cite{NPS}.

The aim of this note is to interpret (\ref{eq.conserve}) as a {\it linear} conservation equation solved
 in the duality sense \cite{BJ98}, the product being defined 
afterwards, following the strategy introduced in \cite{BJ99} for pressureless gases. Therefore we recall in the 
next section the notion of duality solutions and some useful results.
In Section 3 we state and prove the main result concerning
existence and uniqueness of duality solutions of system (\ref{eq:nl}).
Section 4 is devoted to some examples of applications of this result.

\section{Duality solutions for linear equations}
The notion of duality solutions was introduced in \cite{BJ98} to
give a sense to linear conservation equations
\beq\label{eq.conserve}
\p_t \rho + \p_x(b \rho )=0,
\eeq
when the coefficient $b$ can be discontinuous but satisfies
the so-called one-sided Lipschitz (OSL) condition
\begin{equation}\label{OSLC}
\partial_x b(t,.)\leq \beta(t)\qquad\mbox{for $\beta\in L^1(0,T)$ in the distribution sense}.
\eeq
Duality solutions are defined as weak solutions, the test 
functions being specific Lipschitz solutions to the backward linear transport equation
\begin{equation}
\partial_t p + b(t,x) \partial_x p = 0, 
\quad p(T,.) = p^T \in {\rm Lip}(\RR).
\label{transp}\end{equation}
\begin{definition}
\begin{enumerate}
\item We say that a Lipschitz solution $p$ to (\ref{transp}) is a {\bf reversible solution} 
if $p$ is locally constant on the set 
	$$
{\mathcal V}_e=\Big\{(t,x)\in [0,T] \times \RR;
\ \exists\ p_e\in{\mathcal E},\ p_e(t,x)\not=0\Big\}.
	$$
\item We say that 
$\rho\in C([0,T];{\mathcal M}_{loc}(\RR)-\sigma({\mathcal M}_{loc},C_c))$
is a {\bf duality solution} to (\ref{eq.conserve}) 
if for any $0<\tau\le T$, and any {\bf reversible} solution $p$ to (\ref{transp})
with compact support in $x$,
the function $t\mapsto\int_{\RR}p(t,x)\rho(t,dx)$ is constant on
$[0,\tau]$.
\end{enumerate}
\end{definition}

The most important facts for our purpose concerning duality solutions are gathered in the following theorem.
\begin{theorem}(Bouchut, James \cite{BJ98})\label{ExistDuality}
\begin{enumerate}
\item Given $\rho^\circ \in {\mathcal M}_{loc}(\RR)$, under the 
assumptions (\ref{OSLC}), there exists a unique 
$\rho \in C([0,+\infty[,{\mathcal M}_{loc}(\RR))$,
duality solution to (\ref{eq.conserve}), such that $\rho(0,.)=\rho^\circ$. \\
Moreover, if $\rho^\circ$ is nonnegative, then $\rho(t,\cdot)$ is nonnegative
for a.e. $t\geq 0$. And we have the mass conservation 
$|\rho(t,\cdot)|(\RR) = |\rho^\circ|(\RR)$, for a.e. $t\in ]0,T[$.
\item Backward flow and push-forward: the duality solution satisfies
\beq\label{flow}
\forall\, t\in [0,T], \forall\, \phi\in C_c(\RR),\quad
\int_\RR \phi(x)\rho(t,dx) = \int_\RR \phi(X(t,0,x)) \rho^0(dx),
\eeq
where the {\bf backward flow} $X$ is defined as the unique reversible
solution to
	$$
\p_tX+b(t,x) \p_xX = 0 \quad \mbox{ in } ]0,s[\times\RR, \qquad
X(s,s,x)=x.
	$$
\item There exists a bounded Borel function $\bchapo$, called {\bf
universal representative} of $b$, such that $\bchapo = b$
almost everywhere, and for any duality solution $\rho$,
\beq\label{eq:distrib}
\partial_t \rho + \partial_x(\bchapo\rho) = 0 \qquad \hbox{in the distributional sense.}
\eeq
\end{enumerate}
\end{theorem}
\begin{remark}\label{dual.trans}
A similar notion of duality solution for the transport equation is available
 $\p_t u + b\p_x u = 0$,
and $\rho$ is a duality solution of (\ref{eq.conserve}) iff $u=\int^x\rho$
is a duality solution to transport equation (see \cite{BJ98}).
\end{remark}

We shall need also the following result whose proof can be found in \cite{BJ99} 
(Theorems 3.1 and 3.2)
\begin{theorem}\label{th.entrop}
Let $f\in C^1(\RR)$.
Let $M$ be an entropy solution to the conservation equation
$$
\p_t M + \p_x f(M) = 0,
$$
with nondecreasing initial datum $M^0$. 
Then $\rho:=\p_xM$ is a duality solution to 
$$
\p_t\rho + \p_x(b\rho)=0
$$ 
where we can choose $b=f'(M)$ a.e.
Moreover, for all $t\in ]0,+\infty[$, $\p_xb\leq 1/t$ and its 
universal representative $\bchapo$ satisfies
$
\p_xf(M) = \bchapo \p_x M.
$
\end{theorem}

\section{Duality solutions for weakly nonlinear equations}

We introduce the following notion of duality solution for the coupled system 
(\ref{eq:nl}), inspired by the strategy used in \cite{BJ99} for pressureless gases 
(see also section \ref{VPFP} below).
\begin{definition}\label{dual.inter}
We say that $(\rho,u)$ is a duality solution 
of (\ref{eq:nl}) on $]0,T[$ if there exists a bounded Borel function
$b$ with $\p_xb\leq \alpha \in L^1_{loc}(0,T)$ such that
\begin{enumerate}
\item for all $0<t_1<t_2<T$, $\p_t \rho + \p_x(b\rho)=0$ in the sense of 
duality on $]t_1,t_2[$,
\item we have $\p_{x}u = \rho$ in the weak sense,
\item $b=a(u)$ almost everywhere.
\end{enumerate}
\end{definition}

We underline at once the fact that this definition does not lead to uniqueness as it stands.
Indeed, assume that 
$a$ is a non-increasing $C^1$ function on $\RR$ and take for initial data
$\rho^0=\delta_{x_0}$, a Dirac measure in $x_0\in \RR$.
Looking for a solution as a Dirac mass $\delta_{x_1(t)}$, 
thanks to Remark \ref{dual.trans} we solve the transport equation
with coefficient $a\big(H(x-x_1(t))\big)$, where $H$ denotes the Heaviside function.
Then $\big(\delta_{x_1(t)},H(x-x_1(t))\big)$ is a duality solution of (\ref{eq:nl}) 
in the sense of Definition \ref{dual.inter}, provided that $x_1(0)=x_0$ 
and that the admissibility condition $a(1)<x_1'(t)<a(0)$ holds. 
Thus we have an infinite family of duality solutions.

Therefore the main result of this note is to explain how a more precise description of the product $b\rho$ in the scalar 
conservation equation allows to recover uniqueness. 
It is actually given in a very naive way by
writing $a(u)\rho=a(u)\partial_xu=\partial_xA(u)$, where $A$ is an antiderivative of $a$. 
This choice can be justified in a more rigourous way when the system (\ref{eq:nl}) is obtained as the hydrodynamic
limit of a kinetic system, as it is the case both in \cite{jv} and \cite{NPS}. 
It turns out that the previous formal computation is correct at the kinetic level, so that the flux
$J_\eps:=\int\xi f_\eps(\xi)\,d\xi$, where $f_\eps$ is the distribution function of particles, actually converges 
to $J=\partial_xA(u)$, which defines the flux of the conservation equation.
The point now is to justify that this can be used to solve the conservation equation in the duality sense.

\begin{theorem}\label{th1}
Let $\rho^0\in {\mathcal M}_{loc}(\RR)$, $\rho^0\geq 0$.
There exists a unique duality solution $(\rho,u)$ to
the non local interaction equation (\ref{eq:nl})
in the sense of Definition \ref{dual.inter}, which satisfies 
$\bchapo \rho = \p_x(A(u))$ where $A$ is an antiderivative of $a$.

Moreover, if $a$ is a non-increasing function (i.e. in the attractive case), 
there exists a flow associated to $a(u)$, in other words, 
there exists a Lipschitz function $X$ such that $\rho(t)=X(t)_\#\rho^0$.
\end{theorem}

\smallskip

\noindent{\bf Proof.} 
Let us denote by $u$ the entropy solution of equation 
\beq\label{eq:entrop}
\p_t u + \p_x A(u) = 0,
\eeq
with initial data $u^0:=\int^x\rho^0(dx)$. 
From Theorem  \ref{th.entrop}, $\rho:=\p_xu$ is a duality solution of 
$\p_t\rho+\p_x(b\rho)=0$ where we can choose $b=a(u)$ almost everywhere and it
is the unique solution satisfying $\bchapo \rho = \p_x A(u)$. 
Indeed if we have two such nonnegative solutions $\rho_1=\p_xu_1$ 
and $\rho_2=\p_xu_2$, then $u_1$ and $u_2$ are monotonous solutions 
of (\ref{eq:entrop}). Thus they are entropy solutions of this 
scalar equation and $u_1=u_2$.

In the attractive case, the entropy solution $u$ is nondecreasing. 
Therefore $\rho:=\p_xu$ is nonnegative and 
$\p_xa(u)=a'(u)\p_xu\leq 0$ since in the attractive case $a$ is 
non-increasing.
Thus the velocity field $a(u)$ satisfies the OSL condition (\ref{OSLC})
and from Theorem \ref{ExistDuality} $(ii)$ there exist a backward flow 
$X$ such that (\ref{flow}) is satisfied.

In the general case, we can apply the classical Oleinik entropy 
condition and get that $\p_xb\leq 1/t$.
Then the solution is defined on all $]t_1,t_2[$ for $0<t_1<t_2<T$ 
and the flow cannot be defined up to $0$.

\section{Examples}

\subsection{Positive chemotaxis}
Equation (\ref{eq:nl}) for a non-increasing function $a$ 
can be obtained from a hydrodynamical limit of a kinetic model 
describing positive chemotaxis (see e.g. \cite{dolschmeis,jv}).
Thus from Theorem \ref{th1}, there exists a flow $X$ such that
$\rho=X_\#\rho^0$.
Let us first come back to the example in 
subsection 3.1: we assume that $a$ is a given non-increasing 
$C^1$ function and take $\rho^0=\delta_{x_0}$. 
Then we solve the Riemann problem 
$$
\p_t u + \p_x A(u) = 0, \qquad u(t=0,x)=H(x-x_0),
$$
where $A$ is a concave function. Then the entropy solution is 
given by $u(t,x)=H(x-x_1(t))$ where the Rankine-Hugoniot 
condition implies $x_1'(t)=A(1)-A(0)$.
Thus the unique duality solution in the sense of Theorem \ref{th1}
is given by $(\rho,u)=(\delta_{x_1(t)},H(x-x_1(t)))$ where 
$x_1(t)=x_0+(A(1)-A(0))t$.

On the other hand, if we look for a solution in the form $\rho(t)=\delta_{x_1(t)}$,
then $u(t,x)=H(x-x_1(t))$. Integrating equation (\ref{eq:distrib})
we get that $\bchapo \rho=-\p_t u=x_1'(t)\delta_{x_1(t)}$.
By deriving in the distribution sense $A(u)$, we get that
the definition of the product in Theorem \ref{th1} $\bchapo \rho=\p_xA(u)$ 
is satisfied if and only if $x_1'(t)=A(1)-A(0)$, thus we recover the Rankine-Hugoniot condition. 
Hence the definition of the product allows to select one solution among those found in subsection 3.1.
It gives more generally the dynamic of aggregates, which are modelled by a sum of Dirac masses
$\rho^0=\sum m_i\delta_{x_i}$. A similar computation gives the velocity $x'_i(t)=\big(A(\sum^im_j)-A(\sum^{i-1}m_j)\big)/m_i$.
Notice that the velocity of each aggregate is defined by a local equation, despite the fact that the initial equation is non local.
In the particular case where $A$ is strictly concave, aggregates collapse 
in finite time.

\subsection{High field limit of Vlasov--Poisson--Fokker--Planck}\label{VPFP}
In \cite{NPS}, the authors prove that solutions to the 
Vlasov--Poisson--Fokker--Planck system converge in the high field limit
to solutions of (\ref{eq:nl}) where $a(u)=u$ in the repulsive case 
and $a(u)=-u$ in the attractive case. 
To do so, the authors define a weak product $\rho u$, 
which can be proved to coincide with the one used here. 
Applying Theorem \ref{th1} we can recover the result stated in Theorem 2 of \cite{NPS}~:
there exists a unique global in time solution of (\ref{eq:nl}) 
in the distribution sense such that the product $\rho u = \pm u^2/2 $. 
Moreover, in the attractive case, there exists a flow $X$ such that $\rho=X_\#\rho^0$, and 
the dynamics of aggregates is similar to the one of chemotaxis.
In the general case, the Oleinik entropy condition gives that $\rho\leq 1/t$.
Finally, we notice that the result of \cite{NPS} has been extended in 
two dimensions by Poupaud in \cite{poupaud} by using defect measures 
to define the product of $\rho$ by $u$. However, there is no uniqueness 
of solutions.

To conclude, we focus on the connection between pressureless gases and the Vlasov--Poisson--Fokker--Planck limit
which is mentioned in \cite{NPS}. The pressureless gases system reads
\beq\label{pressureless}
\p_t \rho + \p_x (\rho v) = 0, \qquad
\p_t (\rho v) + \p_x (\rho v^2) = 0.
\eeq
Bouchut and James in \cite{BJ99} introduced the notion of duality solution to (\ref{pressureless}):
\begin{definition}\label{p0.def}
We say that a couple $(\rho,q)$, $\rho,q \in C([0,T[;{\mathcal M}_{loc}(\RR))$,
$\rho\geq 0$, is a duality solution to (\ref{pressureless}) if there exists 
$b\in L^\infty(]0,T[\times \RR)$ and $\alpha\in L^1_{loc}(]0,T[)$
satisfying $\p_xb\leq \alpha$ in $]0,T[\times \RR$ such that
\begin{enumerate}
\item For all $0<t_1<t_2<T$, we have in the sense of duality on 
 $]t_1,t_2[\times \RR$
$$
\p_t \rho +\p_x(b\rho) = 0, \qquad \p_tq + \p_x(bq)=0;
$$
\item $\bchapo \rho = q$.
\end{enumerate}
\end{definition}
The existence result for the Cauchy problem with initial data $(\rho^0,q^0)$
 strongly exploits the relationships between (\ref{pressureless}) and the conservation law
$\p_t u + \p_xA(u)=0$, where $\rho=\p_x u$, $q=\p_xA(u)$ and $A$ is determined by $\rho$ and $q$.
Uniqueness follows if $A$ can be defined by the initial data, which enforces additional conditions on $(\rho^0,q^0)$.

In the context of Vlasov--Poisson--Fokker--Planck, the function $A$ is given: $A(u)=-u^2/2$ (attractive case) or $A(u)=u^2/2$ (repulsive case).
Therefore we propose the following variant to the results of \cite{BJ99}.
\begin{theorem}\label{th.p0}
Let $\rho^0\in {\mathcal M}_{loc}(\RR)$, $\rho^0\geq 0$ and 
$A\in C^1(\RR)$. Define $q^0=\p_x A(u^0)$, where  $u^0=\int^x \rho^0$. 
Then there exists a duality solution $(\rho^0,q^0)$ of the pressureless
gases system (\ref{pressureless}) in the sense of Definition \ref{p0.def}.
Moreover this solution is the unique duality solution 
which satisfy the relation
$\bchapo \rho = \p_x [A(u)]$, where $u = \int^x \rho(dx)$. 
\end{theorem}
This theorem is proved in the same way as Theorem \ref{th1}. 
The solution to the high-field limit of the Vlasov--Poisson--Fokker--Planck system
obtained in \cite{NPS} is therefore the unique duality solution to (\ref{pressureless}) given by Theorem \ref{th.p0} for the corresponding
$A$.



\end{document}